\documentclass[12pt,twoside]{article}

\usepackage{graphicx,graphics, amsthm, amsfonts, amsbsy,amssymb, amsmath, cite, array,arydshln, hyperref, float,tikz}

\everymath{\displaystyle}

\usepackage[margin=1in]{geometry}

\allowdisplaybreaks

\newtheorem{theorem}{Theorem}[section]
\newtheorem{remark}[theorem]{Remark}

\newtheorem{lemma}[theorem]{Lemma}
\newtheorem{corollary}[theorem]{Corollary}

\newtheorem{proposition}[theorem]{Proposition}

\newtheorem{conjecture}{Conjecture}

\def\qed{\nolinebreak\hfill\rule{.2cm}{.2cm}\par\addvspace{.5cm}}

\begin{document}
\title{A note on the eigenvalues of zero divisor graphs associated with commutative rings}
\author{ Bilal Ahmad Rather$^{a,}$\footnote{Orcid: 0000-0003-1381-0291}\\ 
	$^{a}${\it Department of Mathematical Sciences, College of Science}, \\
	{\it United Arab Emirate University, Al Ain 15551, Abu Dhabi, UAE}\\
	$ ^{a} $\texttt{bilalahmadrr@gmail.com} 
}
\date{}

\pagestyle{myheadings} \markboth{Bilal Ahmad Rather}{A note on the eigenvalues of zero divisor graphs associated with commutative rings}
\maketitle
\vskip 5mm
\noindent{\footnotesize \textbf{ Abstract}.} For a commutative ring $R,$ with non-zero zero divisors $Z^{\ast}(R)$. The zero divisor graph $\Gamma(R)$ is a simple graph with vertex set $Z^{\ast}(R)$, and two distinct vertices $x,y\in V(\Gamma(R))$ are adjacent if and only if $x\cdot y=0.$ In this note, we provide counter examples to the eigenvalues, the energy and the second Zagreb index related to zero divisor graphs of rings obtained in [Johnson and Sankar, J. Appl. Math. Comp. (2023), \cite{johnson}]. We correct  the eigenvalues (energy) and the Zagreb index result  for the zero divisor graphs of ring $\mathbb{Z}_{p}[x]/\langle x^{4} \rangle.$ We show that for any prime $p$, $\Gamma(\mathbb{Z}_{p}[x]/\langle x^{4} \rangle)$ is non-hyperenergetic and for  prime $p\geq 3$, $\Gamma(\mathbb{Z}_{p}[x]/\langle x^{4} \rangle)$ is hypoenergetic. We give a formulae for the topological indices of $\Gamma(\mathbb{Z}_{p}[x]/\langle x^{4} \rangle)$ and show that its Zagreb indices satisfy Hansen and Vuki$\check{c}$cevi\'c conjecture \cite{hansen}.
\vskip 3mm



\section{Introduction}
\paragraph{}
A finite simple graph $G=G(V(G),E(G))$ consists of a vertex set $V(G)$ and an edge set $ E(G) $. The number $|V(G)|=n$ is \textit{order} and the number  $|E(G)|=m$ is \textit{size} of $G$. By $u\sim v$, we denote $ u $ and $ v $ are adjacent. The set of vertices $N_{G}(v)$ adjacent to $v$ is the \textit{neighbourhood} (open neighbourhood) of $v.$ The number $|N_{G}(v_{i})| $ is \textit{degree} $ d_{i} $ (or $ d_{v_{i}} $) of $ v_{i} $. 
An \emph{independent} set is a collection of non-adjacent vertices and a \emph{clique} is a subset of $V$ such that any two vertices in it are adjacent. 
A graph $ G $ is regular if $ d_{i} $ is same for every $i.$  The \textit{adjacency matrix} $A(G)$ is matrix of $0$'s and $1$'s such that its $(i,j)$-th entry is one if $v_{i}\sim v_{j}$, and $0$, otherwise. $A(G)$  is a real and symmetric matrix, we order its eigenvalues as 
\[ \lambda_{n}(G)\geq \lambda_{n-1}(G)\geq \dots\geq \lambda_{2}(G)\leq \lambda_{1}(G), \] and the multiset of $ \lambda_{i}(G) $'s is its \emph{spectrum}.  
The \emph{energy} of $ G $ is defined as
\[ \mathcal{E}(G)=\sum_{i=1}^{n}|\lambda_{i}(G)|. \]
The energy $ \mathcal{E}(G) $ of a graph $G$ is very well studied and was defined by Gutman in 1978. The characteristic polynomial of a matrix $M$ is $\Phi(M,\lambda)=\det(\lambda I-M),$ where $\det$ is determinant of $M$ and $I$ is the identity matrix.
However the intension of energy was already there, since in 1930s,  Huckel's  Molecular Orbital Theory allows us to predict the approximate values of energies associated with orbitals of electrons in hydrocarbons. The method works with the process that the Hamiltonian operator can be written as a linear combination of certain orbitals, and it provides us required energies with the help of the time-independent Schrödinger equation. But Gunthard and Primes discovered in 1956  that the matric method used by Huckel  is a first degree characteristic polynomial of $A(G)$ for a  graph structure $G$ associated to molecule. There is a vast literature about $A(G)$ and its energy, see \cite{shi, filipovski, gutman border, gutman2010}.

For a commutative ring $ R $ with non-zero identity. An element $ 0\neq x\in R $ is known as the zero divisor of $ R $ if we can find $ 0\neq y\in R $ such that $ x\cdot y=0. $ Beck \cite{ib} introduced the concept of the zero divisor graphs of commutative rings and included $ 0 $ in the definition. While Anderson and Livingston \cite{al} excluded $0$ and  altered the definition of zero divisor graphs, they defined edges between two non zero zero divisors $x$ and $y$ if and only if $x\cdot y=0$.  The adjacency eigenvalues of zero divisor graphs were first discussed by Young \cite{my}. The spectral theory was further elaborated in \cite{magi1, monius}. Singh and Bhat \cite{singh, singh1} gave several topological indices and the energy (Laplacian energy) of zero divisor graphs (and their line graphs) of commutative rings. 
Motivated by the work of authors \cite{singh, singh1, johnson}, we further investigate the zero divisor graphs and correct some already published results in \cite{johnson}.

We use standard notations, complete graph is denoted by $ K_{n}$, and notations/terminology, the readers are referred to \cite{cds,hj,roman}.

The rest of the paper is organized as follows. In Section \ref{section 2}, we give the eigenvalues, the eigenvectors and the energy of the zero divisor graphs of $\mathbb{Z}_{p}[x]/\langle x^{4} \rangle,$ thereby correcting the earlier published results by Johnson and  Sankar in \cite{johnson}. We prove that for any prime $p$, $\Gamma(\mathbb{Z}_{p}[x]/\langle x^{4} \rangle)$ is non-hyperenergetic and for  prime $p\geq 3$, $\Gamma(\mathbb{Z}_{p}[x]/\langle x^{4} \rangle)$ is hypoenergetic. In Section \ref{section 3}, we give the topological indices of $\Gamma(\mathbb{Z}_{p}[x]/\langle x^{4} \rangle),$ which in turn corrects the result related to second Zagreb index recently published in \cite{johnson}. We also show that the Zagreb indices of $\Gamma(\mathbb{Z}_{p}[x]/\langle x^{4} \rangle)$ satisfy the
Hansen and Vuki$\check{c}$cevi\'c conjecture. We end up the article with a conclusion.

\section{Eigenvalues and energy of zero divisor graphs} \label{section 2}
\paragraph{} 
For the sake of completeness and the typos in the construction of $\Gamma(\mathbb{Z}_{p}[x]/\langle x^{4} \rangle), $ we discuss its structural process. 
For a prime $ p $, the structure of zero divisor graph $G\cong \Gamma(\mathbb{Z}_{p}[x]/\langle x^{4} \rangle)$ (see, \cite{johnson}) of $\mathbb{Z}_{p}[x]/\langle x^{4} \rangle$ can be constructed as given below:

We partition the vertex set $V(G)$ of $G$ into following sets
\begin{align*}
	A&=\Big\{ kx^{3} | k=1,2,\dots, (p-1) ~\text{and}~ p\nmid k\Big\},\\
	B&=\Big\{ lx^{3}+lx^{2} | l=1,2,\dots, (p-1) ~\text{and}~ p\nmid l \Big\},\\
	C&=\Big\{ mx^{3}+lx^{2}+kx | m=1,2,\dots, (p-1) ~\text{and}~ p\nmid m \Big\}.
\end{align*} 
The elements of $A$ are $\Big\{  x^{3}, 2x^{3}, 3x^{3}, \dots, (p-2)x^{3}, (p-1)x^{3} \Big\}$ and its cardinality is $|A|=p-1.$ Similarly, $|B|=p(p-1)$ and that of $C$ is $p(p^{2}-p).$ Further, for arbitrary $k_{i}x^{3}$ and $k_{j}x^{3}$ in $A,$ we have
\[ (k_{i}x^{3})\cdot (k_{j}x^{3})=k_{i}k_{j}x^{6}\equiv 0 \mod x^{4}. \]
Thus, it follows that each element of $A$ is adjacent to the remaining elements contained inside $A.$ Again for any two elements say $l_{i}x^{3}+k_{j}x^{2}$ and $l_{j}x^{3}+k_{i}x^{2}$ in $B,$ we have
\[ (l_{i}x^{3}+k_{j}x^{2})\cdot l_{j}x^{3}+k_{i}x^{2}=l_{i}l_{j}x^{6}+(l_{i}k_{i}+l_{j}k_{j})x^{5}+(k_{i}k_{j})x^{4} \equiv 0 \mod x^{4}. \]
It implies that each vertex in $B$ is adjacent to every other vertex of $B,$ that is, the vertices of $B$ form a clique of size $p^{2}-p$. Again for $m_{i} x^{3} +l_{i} x^{2} +k_{i} x \in C$ and $m_{j} x^{3} +l_{j} x^{2} +k_{j} x \in C,$ then we have
\begin{align*}
	 (m_{i} x^{3}&+l_{i} x^{2} +k_{i} x)\cdot (m_{j} x^{3} +l_{j} x^{2} +k_{j} x)\\
	 &= m_{i}m_{j}x^{6} + (l_{i}m_{j} +l_{j}m_{i})x^{5} +(k_{i}m_{j}+l_{i}l_{j}+m_{i}k_{j})x^{4}+(k_{i}l_{j} +l_{i}k_{j})x^{3}+k_{i}k_{j}x^{2}\\
	 &\not \equiv 0 \mod x^{4}.
\end{align*}
Thus, the vertices in $C$ are not adjacent to themselves, so they form an independent set of cardinality $p^{3}-p^{2}.$ Again, for $k_{i} x^{3} \in A$ and $l_{i}x^{3} + k_{j} x^{2}\in B,$ we see that 
\[ (k_{i} x^{3})\cdot(l_{i} x^{3} + k_{j} x^{2}) = l_{i} k_{i} x^{6} + k_{i} k_{j} x^{5} \equiv 0 \mod x^{4}, \] 
so it implies that each vertex of $A$ is adjacent to every vertex in $B$. Likewise for any $k_{i} x^{3}\in A$ and $m_{i} x_{3} + l_{i} x_{2} + k_{j} x\in C,$ we see that 
$$ (k_{i} x^{3})\cdot (m_{i} x_{3} + l_{i} x_{2} + k_{j} x)=m_{i} k_{i} x^{6}+l_{i} k_{i} x^{5}+k_{i} k_{ j}x^{4} \equiv 0 \mod x^{4}.$$ So each vertex in $A$ is connected to every vertex in $C.$ Lastly, if $l_{i} x^{3} + k_{i} x^{2} \in B$ and $m_{j} x_{3} + l_{j} x_{2} + k_{j} x \in C,$ then 
$$ (l_{i} x^{3} + k_{i} x^{2})\cdot (m_{j} x_{3} + l_{j} x_{2} + k_{j} x )=l_{i}m_{j} x^{6} + (k_{i}m_{j} +l_{i} l_{j})x^{5} + (k_{i} l_{j} +l_{i} k_{j})x^{4} +k_{i}k^{j}x^{3}\not\equiv 0 \mod x^{4}.$$ Thus, t follows that the vertices in $B$ are not adjacent to any vertex in $C$ (while authors in \cite{johnson} have written opposite of it, may be a typo). This gives us the structure of $G$ completely. We illustrate this process with the help of an example (same as given in \cite{johnson}). 

\begin{figure}[H]
	\centering
	\scalebox{.3}{\includegraphics{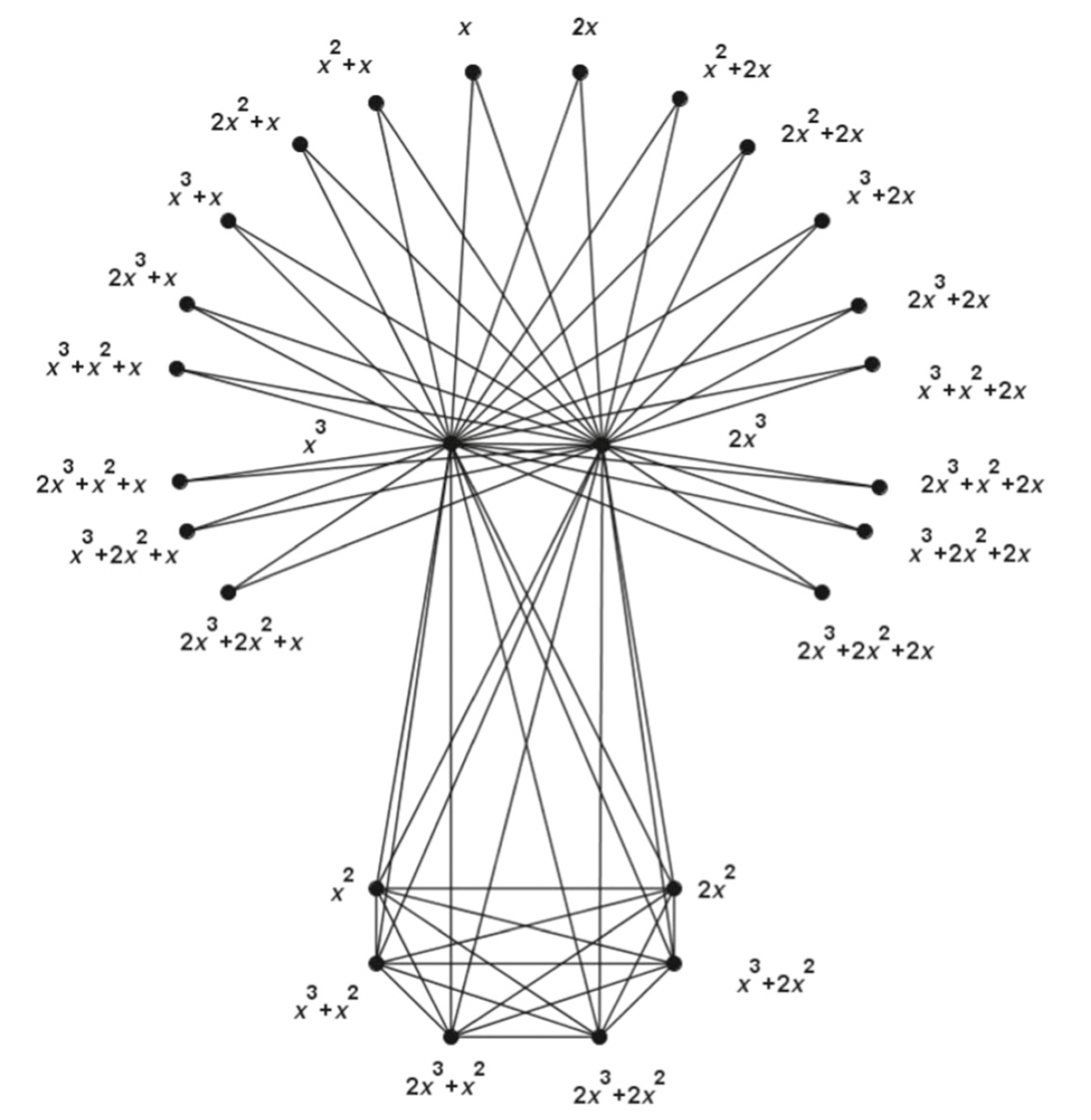}}
	\caption{Zero divisor graph of $\Gamma \big(\mathbb{Z}_{3}[x]/\langle x^{4} \rangle \big)$. }
	\label{fig 1}
\end{figure}

For $R\cong \mathbb{Z}_{3}[x]/\langle x^{4} \rangle,$ with $p=3.$ the vertex set of $G \cong \Gamma \big(\mathbb{Z}_{3}[x]/\langle x^{4} \rangle \big)$ are given as
\begin{equation}\label{sizes of ABC for example}
	\begin{aligned}
	A &= \{x^{3}, 2x^{3}\},\\
	B &= \{x^{2}, 2x^{2}, x^{3} + x^{2}, x^{3} + 2x^{2}, 2x^{3} + x^{2}, 2x^{3} + 2x^{2}\},\\
	C &= \{x, 2x, x^{2} + x, x^{2} + 2x, 2x^{2} + x, 2x^{2} + 2x, x^{3} + x, x^{3} + 2x, 2x^{3} + x,2x^{3} + 2x,\\
	& x^{3} + x^{2} + x, x^{3} + x^{2} + 2x, x^{3} + 2x^{2} + x, x^{3} + 2x^{2} + 2x,2x^{3} + x^{2} + x, 2x^{3} + x^{2} + 2x,\\
	& 2x^{3} + 2x^{2} + x, 2x^{3} + 2x^{2} + 2x\}.
\end{aligned}
\end{equation}
The pictorial representation of $G$ is shown in Figure \ref{fig 1}, where degree of each vertex in $C$ is $2,$ degree of each vertex in $B$ is $7$ and the two vertices of $A$ are of dominating vertices (adjacent to all vertices) of degree $ 25. $

\begin{remark}
	From the structure of $G\cong \Gamma \big(\mathbb{Z}_{p}[x]/\langle x^{4} \rangle \big),$ it follows that there are $p^{3}-p^{2}$ independent vertices, $p-1$ dominating vertices  and $p^{2}-p$ vertices whose induced subgraph is a clique. Overall, there the clique size of $G$ is $p^{2}-1.$
\end{remark}

Now, we consider the results about the energy of $\Gamma(\mathbb{Z}_{p}[x]/\langle x^{4} \rangle)$ and the adjacency matrix of $\Gamma(\mathbb{Z}_{3}[x]/\langle x^{4} \rangle).$

\begin{theorem}[Theorem 3.2, \cite{johnson}]\label{jhonson energy correction}
	Let $R\cong \mathbb{Z}_{p}[x]/\langle x^{4}\rangle$ be a zero divisor graph with prime $p.$ Then energy of $\Gamma(\mathbb{Z}_{p}[x]/\langle x^{4} \rangle)$ is \[ \mathcal{E}(\Gamma(\mathbb{Z}_{p}[x]/\langle x^{4} \rangle))\geq \frac{1}{2}\Big(p^{5} + p^{4}-4p^{3} + 3p^{2}-5p + 4\Big). \]
\end{theorem}
We outline some lines of the proof here.

\textit{In the proof of the theorem, the characteristic polynomial of $G\cong \Gamma(\mathbb{Z}_{p}[x]/\langle x^{4} \rangle)$ and other steps are written as below:
\begin{equation}\label{chracteristic polynomial jhonson}
	\begin{aligned}
\Phi(g,\lambda)&=\Big(\lambda - p^{3} + 1\Big)^{(p-2)}\Big(\lambda - p + 1\Big)^{p^{2}(p-1)-1}\left (\lambda - \frac{1}{2}\big(p^{4} - 2p^{3} + 3p\big) \right )^{p(p-1)-1}\\
&\times \Big(\lambda^{3} -\frac{1}{2}(p^{4} -2p^{2} + 5p - 2)\lambda^{2} + \frac{1}{2}p(p -1)(p^{2} + p + 1)(p^{3}- 2p^{2}- 2p + 5)\lambda\\
&- \frac{1}{2}p(p-1)^{2}(p^{4}-2p^{3} - 2p^{2} + 3p + 2)\Big).
\end{aligned}
\end{equation}
We have $\alpha + \beta + \gamma = \frac{1}{2}(p^{4}-2p^{2} + 5p- 2),$ the authors did not mention the entities $\alpha, \beta$ and $\gamma,$ but it seems they are the zeros of  of the polynomial $\lambda^{3} -\dfrac{1}{2}(p^{4} -2p^{2} + 5p - 2)\lambda^{2} + \frac{1}{2}p(p -1)(p^{2} + p + 1)(p^{3}- 2p^{2}- 2p + 5)\lambda- \frac{1}{2}p(p-1)^{2}(p^{4}-2p^{3} - 2p^{2} + 3p + 2).$
Next, the energy of $G$ is given as
\begin{align}\nonumber
	\mathcal{E}(\Gamma(\mathbb{Z}_{p}[x]/\langle x^{4} \rangle))&=\sum_{i=1}^{n}|\lambda_{i}|=(p-2)|p^{3}-1|+(p^{2}(p-1)-1)|p-1|\\\nonumber
	&+(p(p-1)-1)\left|\frac{1}{2}\big(p^{4}-2p^{3}+3p\big)\right| +|\alpha|+|\beta|+|\gamma|\\\label{eq 1}
	&\geq \frac{1}{2}\big( p^{5}-4p^{3} + 5p^{2}-10p + 6\big)+\alpha+\beta+\gamma\\\label{eq 2}
	&\geq \frac{1}{2}\big(p^{5} + p^{4}-4p^{3} + 3p^{2}- 5p + 4\big).
\end{align}}
\qed

We will illustrate Theorem \ref{jhonson energy correction} and the characteristic polynomial given above with the help of an example.
Let $G\cong \Gamma(\mathbb{Z}_{3}[x]/\langle x^{4} \rangle)$ be the zero divisor graph as given in Figure \ref{fig 1}. We will first write its adjacency matrix, since the matrix given in \cite{johnson} but that is not correct (may be typo error), calculate its characteristic polynomial and deduce its corresponding energy. We index the vertices first from $A,$ then $B$ and finally $C,$ where the values of $A,B$ and $C$ are given in \eqref{sizes of ABC for example}.
\begin{equation}\label{matrix of example}
	\begin{scriptsize}
			M=\left( \begin{array}{cccccccccccccccccccccccccccc}
		0 & 1 & 1 & 1 & 1 & 1 & 1 & 1 & 1 & 1 & 1 & 1 & 1 & 1 & 1 & 1 & 1 & 1 & 1 & 1 & 1 & 1 & 1 & 1 & 1 & 1 \\
		1 & 0 & 1 & 1 & 1 & 1 & 1 & 1 & 1 & 1 & 1 & 1 & 1 & 1 & 1 & 1 & 1 & 1 & 1 & 1 & 1 & 1 & 1 & 1 & 1 & 1 \\
		1 & 1 & 0 & 1 & 1 & 1 & 1 & 1 & 0 & 0 & 0 & 0 & 0 & 0 & 0 & 0 & 0 & 0 & 0 & 0 & 0 & 0 & 0 & 0 & 0 & 0 \\
		1 & 1 & 1 & 0 & 1 & 1 & 1 & 1 & 0 & 0 & 0 & 0 & 0 & 0 & 0 & 0 & 0 & 0 & 0 & 0 & 0 & 0 & 0 & 0 & 0 & 0 \\
		1 & 1 & 1 & 1 & 0 & 1 & 1 & 1 & 0 & 0 & 0 & 0 & 0 & 0 & 0 & 0 & 0 & 0 & 0 & 0 & 0 & 0 & 0 & 0 & 0 & 0 \\
		1 & 1 & 1 & 1 & 1 & 0 & 1 & 1 & 0 & 0 & 0 & 0 & 0 & 0 & 0 & 0 & 0 & 0 & 0 & 0 & 0 & 0 & 0 & 0 & 0 & 0 \\
		1 & 1 & 1 & 1 & 1 & 1 & 0 & 1 & 0 & 0 & 0 & 0 & 0 & 0 & 0 & 0 & 0 & 0 & 0 & 0 & 0 & 0 & 0 & 0 & 0 & 0 \\
		1 & 1 & 1 & 1 & 1 & 1 & 1 & 0 & 0 & 0 & 0 & 0 & 0 & 0 & 0 & 0 & 0 & 0 & 0 & 0 & 0 & 0 & 0 & 0 & 0 & 0 \\
		1 & 1 & 0 & 0 & 0 & 0 & 0 & 0 & 0 & 0 & 0 & 0 & 0 & 0 & 0 & 0 & 0 & 0 & 0 & 0 & 0 & 0 & 0 & 0 & 0 & 0 \\
		1 & 1 & 0 & 0 & 0 & 0 & 0 & 0 & 0 & 0 & 0 & 0 & 0 & 0 & 0 & 0 & 0 & 0 & 0 & 0 & 0 & 0 & 0 & 0 & 0 & 0 \\
		1 & 1 & 0 & 0 & 0 & 0 & 0 & 0 & 0 & 0 & 0 & 0 & 0 & 0 & 0 & 0 & 0 & 0 & 0 & 0 & 0 & 0 & 0 & 0 & 0 & 0 \\
		1 & 1 & 0 & 0 & 0 & 0 & 0 & 0 & 0 & 0 & 0 & 0 & 0 & 0 & 0 & 0 & 0 & 0 & 0 & 0 & 0 & 0 & 0 & 0 & 0 & 0 \\
		1 & 1 & 0 & 0 & 0 & 0 & 0 & 0 & 0 & 0 & 0 & 0 & 0 & 0 & 0 & 0 & 0 & 0 & 0 & 0 & 0 & 0 & 0 & 0 & 0 & 0 \\
		1 & 1 & 0 & 0 & 0 & 0 & 0 & 0 & 0 & 0 & 0 & 0 & 0 & 0 & 0 & 0 & 0 & 0 & 0 & 0 & 0 & 0 & 0 & 0 & 0 & 0 \\
		1 & 1 & 0 & 0 & 0 & 0 & 0 & 0 & 0 & 0 & 0 & 0 & 0 & 0 & 0 & 0 & 0 & 0 & 0 & 0 & 0 & 0 & 0 & 0 & 0 & 0 \\
		1 & 1 & 0 & 0 & 0 & 0 & 0 & 0 & 0 & 0 & 0 & 0 & 0 & 0 & 0 & 0 & 0 & 0 & 0 & 0 & 0 & 0 & 0 & 0 & 0 & 0 \\
		1 & 1 & 0 & 0 & 0 & 0 & 0 & 0 & 0 & 0 & 0 & 0 & 0 & 0 & 0 & 0 & 0 & 0 & 0 & 0 & 0 & 0 & 0 & 0 & 0 & 0 \\
		1 & 1 & 0 & 0 & 0 & 0 & 0 & 0 & 0 & 0 & 0 & 0 & 0 & 0 & 0 & 0 & 0 & 0 & 0 & 0 & 0 & 0 & 0 & 0 & 0 & 0 \\
		1 & 1 & 0 & 0 & 0 & 0 & 0 & 0 & 0 & 0 & 0 & 0 & 0 & 0 & 0 & 0 & 0 & 0 & 0 & 0 & 0 & 0 & 0 & 0 & 0 & 0 \\
		1 & 1 & 0 & 0 & 0 & 0 & 0 & 0 & 0 & 0 & 0 & 0 & 0 & 0 & 0 & 0 & 0 & 0 & 0 & 0 & 0 & 0 & 0 & 0 & 0 & 0 \\
		1 & 1 & 0 & 0 & 0 & 0 & 0 & 0 & 0 & 0 & 0 & 0 & 0 & 0 & 0 & 0 & 0 & 0 & 0 & 0 & 0 & 0 & 0 & 0 & 0 & 0 \\
		1 & 1 & 0 & 0 & 0 & 0 & 0 & 0 & 0 & 0 & 0 & 0 & 0 & 0 & 0 & 0 & 0 & 0 & 0 & 0 & 0 & 0 & 0 & 0 & 0 & 0 \\
		1 & 1 & 0 & 0 & 0 & 0 & 0 & 0 & 0 & 0 & 0 & 0 & 0 & 0 & 0 & 0 & 0 & 0 & 0 & 0 & 0 & 0 & 0 & 0 & 0 & 0 \\
		1 & 1 & 0 & 0 & 0 & 0 & 0 & 0 & 0 & 0 & 0 & 0 & 0 & 0 & 0 & 0 & 0 & 0 & 0 & 0 & 0 & 0 & 0 & 0 & 0 & 0 \\
		1 & 1 & 0 & 0 & 0 & 0 & 0 & 0 & 0 & 0 & 0 & 0 & 0 & 0 & 0 & 0 & 0 & 0 & 0 & 0 & 0 & 0 & 0 & 0 & 0 & 0 \\
		1 & 1 & 0 & 0 & 0 & 0 & 0 & 0 & 0 & 0 & 0 & 0 & 0 & 0 & 0 & 0 & 0 & 0 & 0 & 0 & 0 & 0 & 0 & 0 & 0 & 0 
	\end{array} \right).
	\end{scriptsize}
\end{equation}
By using a computer software Wolfram Mathematica, the characteristic polynomial of $M$ is
\begin{align*}
	\Phi(M,\lambda)=(\lambda^3-6\lambda^{2}-43\lambda+180)(\lambda+1)^{6}(\lambda)^{17},
\end{align*} 
and its
approximated zeros are $ 8.56673, 3.47673, -6.04346, -1$ and $0.$ Thus, the spectrum of $G$ is
\[ \Big\{8.56673, 3.47673, (-1)^{[6]}, 0^{[17]}, -6.04346\Big\}. \]
Now, the energy of $G$ is
\begin{equation}\label{true energy}
	\mathcal{E}(G)=\sum_{i=1}^{26}|\lambda_{i}|=8.56673+3.47673+6+6.04346=24.087806.
\end{equation}

Next, we calculate the characteristic polynomial and the energy of $G$ with the help of Theorem \ref{jhonson energy correction} (Theorem 3.2, \cite{johnson}).
With $p=3$ in \eqref{chracteristic polynomial jhonson}, the characteristic polynomial of $G$ is 
\begin{equation*}
	\begin{aligned}
		\Phi(G,\lambda)&=\Big(\lambda - 3^{3} + 1\Big)^{(3-2)}\Big(\lambda - 3 + 1\Big)^{3^{2}(3-1)-1}\left (\lambda - \frac{1}{2}\big(3^{4} - 2\cdot3^{3} + 3\cdot3\big) \right )^{3(3-1)-1}\\
		&\times \Big(\lambda^{3} -\dfrac{1}{2}(3^{4} -2\cdot3^{2} + 5\cdot3 - 2)\lambda^{2} + \frac{1}{2}3(3 -1)(3^{2} + 3 + 1)(3^{3}- 2\cdot3^{2}- 2\cdot3 + 5)\lambda\\
		&- \frac{1}{2}3(3-1)^{2}(3^{4}-2\cdot3^{3} - 2\cdot3^{2} + 3\cdot3 + 2)\Big),\\
		&=(\lambda-26)(\lambda-2)^{17}(\lambda-18)^{5}(\lambda^{3}-38\lambda^{2}+312\lambda -120).
	\end{aligned}
\end{equation*}
The zeros of the above polynomial are: $26$ with multiplicity one (simple zero), $ 2 $ with order $17$, $ 18 $ with order $5$ and the simple zeros $ 26.3184, 11.2772 $ and $0.404313$. The spectrum of $G$ is
\begin{equation}\label{spectrum example jhonson}
	\Big\{ 26.3184, 26, 18^{[5]}, 11.2772, 2^{[17]},  0.404313\Big\}.
\end{equation}
The spectrum of $G$ given in \eqref{spectrum example jhonson} is not correct, we recall two facts of the adjacency matrix:
\begin{enumerate}
	\item The spectral radius of $G$ cannot exceed maximum degree, so in this case $\lambda_{1}(G)\leq 26 $, while in \eqref{spectrum example jhonson}, $\lambda_{1}(G)$ exceeds $26.$ 
	\item The trace of the adjacency matrix is zero, while in \eqref{spectrum example jhonson}, there are only positive eigenvalues, so their sum cannot be zero.
\end{enumerate}
 Thus the spectrum obtained by Theorem \ref{jhonson energy correction} for $G\cong \Gamma(\mathbb{Z}_{2}[x]/\langle x^{4} \rangle)$ is not correct. Further, by Theorem \ref{jhonson energy correction}, the energy of $G$ is
\begin{align*}
	 \mathcal{E}(G)&=26.3184+26+18\times 5+11.2772+2\times 17+0.404313\\
&= 187.999913.
\end{align*}
Which is different from the computer calculations  given in \eqref{true energy}, and hence is incorrect.

Besides, by the lower bound for the energy given by Theorem \ref{jhonson energy correction} with $p=3$, we have
\[ \mathcal{E}(G)\geq \frac{1}{2}\big(3^{5} + 3^{4}-4\cdot 3^{3} + 3\cdot 3^{2}- 5\cdot 3 + 4\big)= 116,\] which is very high than the actual energy $ \mathcal{E}(G)=24.087806.$ Therefore, with this counter example, neither the characteristic polynomial is correct nor the energy bound is correct. So, Theorem \ref{jhonson energy correction} (Theorem 3.2, \cite{johnson}) is not valid. In the rest of the section, we give the correct form of the characteristic polynomial  of $\Gamma(\mathbb{Z}_{p}[x]/\langle x^{4} \rangle), $ correct its spectrum and present new lower and upper bounds for the energy for the zero divisor family $\Gamma(\mathbb{Z}_{p}[x]/\langle x^{4} \rangle).$\\

We are now in a position to study the spectrum and the energy of $\Gamma(\mathbb{Z}_{p}[x]/\langle x^{4} \rangle)$ and for that
we need some already known results.

\begin{lemma}[\cite{dasDM}]\label{independence eigenvalue}
Let $ G $ be a graph with independent set of vertices $\{v_{1},v_{2},\dots,v_{k}\} $ having same set of neighbors. Then $ G $
has at least $ k-1 $ eigenvalues equal to the eigenvalues $ 0. $
\end{lemma}

For $\alpha=0$ in Theorem 2.1 of \cite{bilalaims}, we have the next very result.
\begin{lemma}[\cite{bilalaims}]\label{eigenvalue of clique number}
	Suppose $G$ is a graph of order $ n $ with vertices $ \{v_{1}, v_{2},\dots, v_{k}\}, (k< n) $ satisfying $ N(v_{i})\setminus \{v_{i}\}=N(v_{j})\setminus \{v_{j}\} $ for all $ i,j\in \{1,2,\dots,k\} $. Then $ G $ has at least $ k-1 $ eigenvalues equal to the eigenvalue $ -1. $ 
\end{lemma}


The following theorem gives the spectrum of the zero divisor graph of $ \Gamma(\mathbb{Z}_{p}[x]/\langle x^{4} \rangle). $

\begin{theorem} \label{spectra of comaximal graphs}
Let $G\cong \Gamma(\mathbb{Z}_{p}[x]/\langle x^{4} \rangle)$ be the zero divisor graph of $ \mathbb{Z}_{p}[x]/\langle x^{4}\rangle $ with $p$ being prime. Then the spectrum of  $ G $ comprises of the eigenvalue $ -1 $ and $0$ with multiplicities $ p^{2}-3 $ and $ p^{3}-p^{2}-1 $, respectively. The other eigenvalues $G$ are the zeros of the polynomial given below
\begin{equation}\label{poly of cubic eq}
 \lambda^3-\lambda^{2}(p^{2}-3)+\lambda(2 - 2p^{2} + 2p^3 - p^4)+p^6-3p^{5}+2p^{4}+p^{3}-p^{2}.
\end{equation}
\end{theorem}
\noindent\textbf{Proof.} We start labelling of vertices form the vertices in $A,$ then  in $B$ and finally  in $C.$ Under this vertex labelling, the adjacency matrix of $ G $ can be written as 
\begin{equation}\label{adjacency matrix of R}
A(G)=\left( \begin{array}{ c c c c c c}
	J_{p-1}-I_{p-1} & J_{(p-1)\times (p^{2}-p)} & J_{(p-1)\times (p^{3}-p^{2}) }\\
	J_{(p^{2}-p)\times (p-1)} & J_{p^{2}-p}-I_{p^{2}-p} &  \textbf{0}_{(p^{2}-p)\times (p^{3}-p^{2}) }\\
	J_{(p^{3}-p^{2})\times (p-1)} & \textbf{0}_{(p^{3}-p^{2})\times (p^{2}-p)} & \textbf{0}_{(p^{3}-p^{2})}
\end{array}\right),
\end{equation}
where $ J $ is the matrix whose each entry is $1$, $ \textbf{0} $ is the zero matrix, and $ I $ is the identity matrix.  Since, the vertices of $ A $  in $ G $ form a clique of size $p-1$ and each such vertex of $ A $ have the same neighbors satisfying the hypothesis of Lemma \ref{eigenvalue of clique number}. Thus, $ G $ has eigenvalue $ -1 $ with multiplicity $ p-2 $ and its associated eigenvectors are
\begin{align*}
	\Big(-1, 1, \underbrace{0, 0,\dots, 0}_{p-3}, \underbrace{0, 0,\dots, 0}_{p^{3}-p}\Big), &\Big(-1, 0, 1, \underbrace{0, 0,\dots, 0}_{p-4}, \underbrace{0, 0,\dots, 0}_{p^{3}-p^{2}}\Big),\\
	&~\vdots\\
	\Big(-1, \underbrace{0, 0,\dots, 0}_{p-4}, -1, 0, \underbrace{0, 0,\dots, 0}_{p^{3}-p}\Big), &\Big(-1, \underbrace{0, 0,\dots, 0}_{p-3},1, \underbrace{0, 0,\dots, 0}_{p^{3}-p}\Big).
\end{align*}
 Again the vertices of $ B $  form a  clique and each vertex $B$ share the common neighbors in $ A, $ so by Lemma \ref{eigenvalue of clique number}, $-1$ is the eigenvalue of $G$ with multiplicity $p^{2}-p-1$. The corresponding eigenvectors of the eigenvalue $-1$ are
 \begin{align*}
 	\Big( \underbrace{0, 0,\dots, 0}_{p-1}, -1, 1, \underbrace{0, 0,\dots, 0}_{p^{2}-p-2 }, \underbrace{0, 0,\dots, 0}_{p^{3}-p^{2}} \Big), 
 	&\Big( \underbrace{0, 0,\dots, 0}_{p-1}, -1, 0, 1, \underbrace{0, 0,\dots, 0}_{p^{2}-p-3 }, \underbrace{0, 0,\dots, 0}_{p^{3}-p^{2}} \Big),\\
 	&~\vdots\\
 	\Big( \underbrace{0, 0,\dots, 0}_{p-1}, -1, \underbrace{0, 0,\dots, 0}_{p^{2}-p-3 },1,0, \underbrace{0, 0,\dots, 0}_{p^{3}-p^{2}} \Big),
 	&\Big( \underbrace{0, 0,\dots, 0}_{p-1}, -1, \underbrace{0, 0,\dots, 0}_{p^{2}-p-2 }, 1, \underbrace{0, 0,\dots, 0}_{p^{3}-p^{2}} \Big).
 \end{align*}
Thus $-1$ is the eigenvalue of $G$ with multiplicity $p^{2}-p-1+p-2=p^{2}-3$. Further, from the construction of zero divisor, the elements in $C$ form an independent set, where each vertex of $C$ have a same neighbourhood. So, by Lemma \ref{independence eigenvalue}, it
  follows that $ 0 $ is the eigenvalue of $ A(G) $ with multiplicity $ p^{3}-p^{2}-1.  $
 The corresponding eigenvectors of the eigenvalue $ 0 $ are
 \begin{align*}
 	\Big( \underbrace{0, 0,\dots, 0}_{p-1},\underbrace{0, 0,\dots, 0}_{p^{2}-p}, -1, 1, \underbrace{0, 0,\dots, 0}_{p^{3}-p^{2}-2 }\Big),& \Big( \underbrace{0, 0,\dots, 0}_{p-1},\underbrace{0, 0,\dots, 0}_{p^{2}-p}, -1, 0, 1, \underbrace{0, 0,\dots, 0}_{p^{3}-p^{2}-3}\Big),\\
 	&~\vdots\\
 	\Big(\underbrace{0, 0,\dots, 0}_{p-1},\underbrace{0, 0,\dots, 0}_{p^{2}-p}, -1, \underbrace{0, 0,\dots, 0}_{p^{3}-p^{2}-3},1,0\Big), &\Big(\underbrace{0, 0,\dots, 0}_{p-1},\underbrace{0, 0,\dots, 0}_{p^{2}-p}, -1, \underbrace{0, 0,\dots, 0}_{p^{3}-p^{2}-2},1\Big).
 \end{align*}
Thus, with this method, we have obtained $p^{3}-p^{2}-1+p^{2}-3=p^{3}-4$ eigenvalues. Next, we are required to find the remaining three eigenvalues of $A(G).$
Let $ X $ be the eigenvector of $ A(G) $ with $ x_{i}=X(v_{i}) $, for $ i=1,2,3,\dots, p^{3}-1.$ Then, it follows that (see, \cite{cds}) every component of $ X $ that corresponds to every vertex of  $ A $ is equal to $ x_{1} $, components of $ X $ corresponding to the vertices of $ B $ is $ x_{2} $,  and the components of $ X $ that corresponds to the vertices $ C $ is taken as $ x_{3} $.
 Thus, with $X=\Big(\underbrace{x_{1}, x_{1},\dots, x_{1}}_{p-1}, \underbrace{x_{2}, x_{2},\dots, x_{2}}_{p^{2}-p}, \underbrace{x_{3}, x_{3},\dots, x_{3}}_{p^{3}-p^{2}}\Big)$,  the eigenequation $ A(G)X=\lambda X $ gives us
 \begin{align*}
 	\lambda x_{1}&=\underbrace{x_{1}+x_{1}+\dots+x_{1}}_{p-2}+\underbrace{x_{2}+x_{2}+\dots+x_{2}}_{p^{2}-p}+\underbrace{x_{3}+x_{3}+\dots+x_{3}}_{p^{3}-p^{2}}\\
 	\lambda x_{2}&=\underbrace{x_{1}+x_{1}+\dots+x_{1}}_{p-1}+\underbrace{x_{2}+x_{2}+\dots+x_{2}}_{p^{2}-p-1}+0 \cdot x_{3}\\
 	\lambda x_{3}&=\underbrace{x_{1}+x_{1}+\dots+x_{1}}_{p-1}+0 \cdot x_{2}+0 \cdot x_{3}
 \end{align*}
 The coefficient matrix of the right hand side of the above system of equations is 
 \begin{equation}\label{Qmat spectra of comaximal graphs}
 Q= \begin{pmatrix}
 	p-2 & p^{2}-p & p^{3}-p^{2}\\
 	p-1 & p^{2}-p-1 & 0\\
 	p-1 & 0 & 0
 \end{pmatrix}.
 \end{equation}
The characteristic polynomial of the above matrix $Q$ given in \eqref{Qmat spectra of comaximal graphs} is
\[ \lambda^3-\lambda^{2}(p^{2}-3)+\lambda(2 - 2p^{2} + 2p^3 - p^4)+p^6-3p^{5}+2p^{4}+p^{3}-p^{2}. \]
The remaining three eigenvalues of $G$ are the zeros of the above polynomial.
   \qed

\begin{remark}
	We tried by computational software to calculate the eigenvalues of the matrix $Q$, but we  were not successful, since they are not easy to locate. Though by intermediate value theorem, we can locate interval where these eigenvalues lie and is of no use for now, especially in establishing bounds for the energy of $G.$
\end{remark}

A graph $G$ is said to be non-singular if the adjacency matrix $det(A(G))\neq 0.$ A graph is said to be without repeated eigenvalues if all its eigenvalues are simple (that is, with multiplicity one). Very often there are papers published related to properties of non-singular graphs and graphs with simple eigenvalues, see \cite{filipovski}.
The following corollary is immediate from Theorem \ref{spectra of comaximal graphs}.
\begin{corollary}
For prime $p,$ the zero divisor graph $\Gamma(\mathbb{Z}_{p}[x]/\langle x^{4} \rangle) $ have the following properties.
\begin{itemize}
	\item[\bf (i)] $ \Gamma(\mathbb{Z}_{p}[x]/\langle x^{4} \rangle)$ is always singular.
	\item[\bf (ii)] $ \Gamma(\mathbb{Z}_{p}[x]/\langle x^{4} \rangle)$ is never with simple eigenvalues for $p\geq3$.
\end{itemize}
\end{corollary}
\noindent\textbf{Proof.} For $p\geq 2,$ the cardinality of $C$ is at least $4,$ so Lemma \ref{independence eigenvalue} implies that $0$ is always the eigenvalues of $\Gamma(\mathbb{Z}_{p}[x]/\langle x^{4} \rangle). $ For the second case with $p=2$, Lemma \ref{eigenvalue of clique number} implies that $-1$ is simple eigenvalue of $\Gamma(\mathbb{Z}_{p}[x]/\langle x^{4} \rangle)$. But Lemma \ref{independence eigenvalue} implies that $0$ is the eigenvalue of $\Gamma(\mathbb{Z}_{p}[x]/\langle x^{4} \rangle)$ with multiplicity $3.$ For $p\geq 3,$ the eigenvalues $-1$ and $0$ are always repeated zeros of the characteristic polynomial of the adjacency matrix of $\Gamma(\mathbb{Z}_{p}[x]/\langle x^{4} \rangle).$ Thus, for any prime the zero divisor graph $\Gamma(\mathbb{Z}_{p}[x]/\langle x^{4} \rangle)$ is always with repeated eigenvalues. \qed


Let $\xi_{1}\geq \xi_{2} \geq \xi_{3}$ be the zeros of polynomial \eqref{poly of cubic eq} or the eigenvalues of the matrix $Q$ given in \eqref{Qmat spectra of comaximal graphs}. The energy of $G$ equals
\begin{equation}\label{energy equation}
	\mathcal{E}(G)=\sum_{i=1}^{n}|\lambda_{i}|=\sum_{i=1}^{p^{2}-3}|-1|+|\xi_{1}|+|\xi_{2}|+|\xi_{3}|= p^{2}-3+\mathcal{E}(Q),
\end{equation}
where $ \mathcal{E}(Q)$ is the energy of $Q.$ We note that $\xi_{i}$ are actually some $\lambda_{i}$'s.

One trivial lower bound from \eqref{energy equation} is
\[ \mathcal{E}(G)\geq p^{2}-3. \]

For non-negative matrices $\xi_{1}$ is same as $\lambda_{1}$ and is always positive, so from \eqref{energy equation}, we have
\[ \mathcal{E}(G)=p^{2}-3+\xi_{1}+|\xi_{2}|+|\xi_{3}|\geq p^{2}-3+\xi_{1}=p^{2}-3+\lambda_{1}, \]
with equality holding if and only if $\xi_{2}=\xi_{3}=0$ or $\xi_{2}$ and $\xi_{3}$ are symmetric towards origin, that is $\xi_{2}=-\xi_{3}.$
By noting that the spectral $ \lambda_{1}$ of any graph is bounded above by maximum degree $\Delta$, that is, $ \lambda_{1}\leq \Delta,$ with equality holding if and only $G$ is $\Delta$ regular. With this observation, we get
\[ \mathcal{E}(G)\leq p^{2}-3+p^{3}-2+|\xi_{2}|+|\xi_{3}|=p^{3}+p^{2}-5+|\xi_{2}|+|\xi_{3}|, \]
since the maximum degree of  $ \Gamma(\mathbb{Z}_{p}[x]/\langle x^{4} \rangle)$ is $n-1=p^{3}-2.$ Equality occurs if and only if $G$ is regular, and in that case $Q$ is rank one matrix along with $\xi_{2}=\xi_{3}=0.$
Now, we establish the lower and the upper bounds for the energy of $\Gamma(\mathbb{Z}_{p}[x]/\langle x^{4} \rangle).$

Let $ \{q_{1}, q_{2}, q_{3},\dots, q_{s}\} $ be the set of positive real numbers and let $ \Pi_{k} $ be the average of products of $ k $-element subset of $ \{q_{1}, q_{2}, q_{3}, \dots , q_{s}\}, $ that is
\begin{align*}
	\Pi_{1} &=\frac{q_{1} + q_{2} + q_{3} + \dots+ q_{s}}{s},\\
	\Pi_{2} &=\frac{1}{\frac{s(s-1)}{2}}\Big(q_{1}q_{2} + q_{1}q_{3} +\dots + q_{1}q_{s} + q_{2}q_{3} +\dots + q_{s-1}q_{s}\Big),\\
	&~\vdots\\
	\Pi_{s} &= q_{1}q_{2}\dots q_{s}. 
\end{align*}
The following Maclaurin symmetric mean inequality \cite{biler} relates $ \Pi_{i} $'s among themselves.
\begin{equation}\label{maclaurin relations}
	\Pi_{1} \geq \Pi_{2}^{\frac{1}{2}} \geq \Pi_{3}^{\frac{1}{3}} \geq \dots \geq \Pi_{s}^{\frac{1}{s}},
\end{equation}
with equalities holding if and only if $ q_{1}=q_{2}=\dots=q_{s}. $

Using these relations, in the next very result, we establish bounds for the energy of $\Gamma(\mathbb{Z}_{p}[x]/\langle x^{4} \rangle)$.
\begin{theorem}\label{energy theorem}
	Let $\Gamma(\mathbb{Z}_{p}[x]/\langle x^{4} \rangle)$ be the zero divisor graph of $\mathbb{Z}_{p}[x]/\langle x^{4} \rangle.$ Then 
	\[
		 \mathcal{E}(G)\geq p^{2}-3+\sqrt{3p^{4}-4p^{3}-2p^{2}+5+6 |p^2 - p^3 - 2 p^4 + 3 p^5 - p^6|^{\frac{2}{3}}}
		  \] and \[
	 \mathcal{E}(G)\leq p^{2}-3+\sqrt{9p^{4}-12p^{3}-6p^{2}+15},
\]
with equalities holding if and only if $ |\xi_{1}(Q)|=|\xi_{2}(Q)|=|\xi_{3}(Q)|.$
\end{theorem}
\noindent\textbf{Proof.} By \eqref{energy equation}, we have 
\begin{equation}\label{energy thrm eq 1}
	 \mathcal{E}(G)=\sum_{i=1}^{n}|\lambda_{i}|= p^{2}-3+\mathcal{E}(Q).
\end{equation}
We will establish bounds for the entity $\mathcal{E}(Q).$ By applying Maclaurin \eqref{maclaurin relations} to the set 
$$\big\{|\xi_{1}(Q)|, |\xi_{2}(Q)|, |\xi_{3}(Q)|\big\},$$ we have
\begin{footnotesize}
	\begin{equation}\label{lemma appl 1}
	\left(\frac{|\xi_{1}(Q)|+|\xi_{2}(Q)|+|\xi_{3}(Q)|}{3}\right)^{2}\geq \frac{1}{\frac{3(3-1)}{2}} \Big( |\xi_{1}(Q)| |\xi_{2}(Q)|+|\xi_{1}(Q)| |\xi_{3}(Q)|+|\xi_{2}(Q)||\xi_{3}(Q)| \Big),
\end{equation}
\end{footnotesize}
that is, 
\begin{align*}
	\Big(|\xi_{1}(Q)|+|\xi_{2}(Q)|+|\xi_{3}(Q)|\Big)^{2} &\geq 3 \Big( |\xi_{1}(Q)| |\xi_{2}(Q)|+|\xi_{1}(Q)| |\xi_{3}(Q)|+|\xi_{2}(Q)||\xi_{3}(Q)| \Big)\\
	&=\frac{3}{2}\left(\Big(\sum_{i=1}^{3}|\xi_{i}(Q)|\Big)^{2}-\sum_{i=1}^{3}\xi_{i}^{2}(Q)\right), 
\end{align*}
that is,
\[\Big(\mathcal{E}(Q)\Big)^{2}=\left(\sum_{i=1}^{3} |\xi_{i}(Q)|\right)^{2}\leq 3 \sum_{i=1}^{3}\xi_{i}^{2}=3 \cdot tr(Q^{2}), \]
where $tr(Q^{2})$ is the trace of $Q^{2}.$  From \eqref{Qmat spectra of comaximal graphs}, the trace of $Q^{2}$ is
\begin{align*}
	 tr(Q^{2})&=(p-2)^{2}+(p^{2}-p-1)^{2}+2(p-1)(p^{2}-p)+2(p-1)(p^{3}-p^{2})\\
	 &=3p^{4}-4p^{3}-2p^{2}+5.
\end{align*}

Substituting it in above expression, we obtain
\[ \mathcal{E}(Q)\leq \sqrt{3 \cdot tr(Q^{2})}=\sqrt{3(3p^{4}-4p^{3}-2p^{2}+5)}. \] Thus by \eqref{energy thrm eq 1}, we get
\[ \mathcal{E}(G)= p^{2}-3+\mathcal{E}(Q)\leq p^{2}-3+\sqrt{9p^{4}-12p^{3}-6p^{2}+15}. \]
Equality holds if and only if equality holds in \eqref{lemma appl 1}, that is, $|\xi_{1}(Q)|=|\xi_{2}(Q)|=|\xi_{3}(Q)|.$
Again from the second inequality of  \eqref{maclaurin relations}, we have
\[\frac{1}{\frac{3(3-1)}{2}} \Big( |\xi_{1}(Q)| |\xi_{2}(Q)|+|\xi_{1}(Q)| |\xi_{3}(Q)|+|\xi_{2}(Q)||\xi_{3}(Q)| \Big)\geq \Big(|\xi_{1}(Q)||\xi_{2}(Q)||\xi_{1}(Q)|\Big)^{\frac{2}{3}},  \]
that is, equivalent to
\[ 2\Big( |\xi_{1}(Q)| |\xi_{2}(Q)|+|\xi_{1}(Q)| |\xi_{3}(Q)|+|\xi_{2}(Q)||\xi_{3}(Q)| \Big)\geq 6\Big(|\xi_{1}(Q)||\xi_{2}(Q)||\xi_{1}(Q)|\Big)^{\frac{2}{3}}, \]
that is, \[ \left(|\xi_{1}(Q)|+|\xi_{2}(Q)|+|\xi_{1}(Q)|\right)^{2}-\Big(|\xi_{1}(Q)|^{2}+|\xi_{2}(Q)|^{2}+|\xi_{3}(Q)|^{2}\Big)\geq  6 |det (Q)|^{\frac{2}{3}}, \]
that is, \[ \mathcal{E}(Q)=\sum_{i=1}^{3}|\xi_{i}(Q)|\geq \sqrt{tr(Q)^{2}+ 6 |det(Q)|^{\frac{2}{3}}}. \]
From \eqref{Qmat spectra of comaximal graphs}, the determinant of $Q$ is $det(Q)=-(p-1)(p^{2}-p-1)(p^{3}-p^{2})=p^2 - p^3 - 2 p^4 + 3 p^5 - p^6.$
Therefore, from \eqref{energy thrm eq 1}, we have
\[ \mathcal{E}(G) \geq p^{2}-3+\sqrt{3p^{4}-4p^{3}-2p^{2}+5+6 |p^2 - p^3 - 2 p^4 + 3 p^5 - p^6|^{\frac{2}{3}}}, \]
with equality holding if and only if $|\xi_{1}(Q)|=|\xi_{2}(Q)|=|\xi_{3}(Q)|.$ \qed
\begin{remark}
	By substituting $p=3,$ in Theorem \ref{energy theorem}, the energy of $G\cong \Gamma(\mathbb{Z}_{3}[x]/\langle x^{4} \rangle)$ satisfies
	\[23.6997\leq  \mathcal{E}(G)\leq 25.1311. \]
	The actual energy of $G$ is $\mathcal{E}(G)=24.0869,$ which is close to the values given by Theorem \ref{energy theorem}.
\end{remark}

Since the energy of $K_{n}$ is $2(n-1)$, and a graph $G$ of order n is said to be hyperenergetic if $\mathcal{E}(G)> 2(n-1)$ and non-hyperenergetic if $E(G) <2(n-1).$ A  graph $G$ with $E(G) = 2(n-1)$ is referred to as borderenergetic. A graph $G$ is called hypoenergetic if $\mathcal{E}(G)<n.$ Over the years many papers were written for borderenergetic, hyperenergetic, hypoenergetic and non-hyperenergetic, for details see \cite{gutman border, hou, ston, gutman2010}. We have the following result in this direction.

\begin{proposition}
	For any prime $p$, the following holds for the zero divisor graph $\Gamma(\mathbb{Z}_{p}[x]/\langle x^{4} \rangle)$ of ring $\mathbb{Z}_{p}[x]/\langle x^{4} \rangle$.
		\begin{itemize}
			\item[\bf (i)] For $p\geq 2$, $\Gamma(\mathbb{Z}_{p}[x]/\langle x^{4} \rangle)$ is non-hyperenergetic.
			\item[\bf(ii)] For $p\geq 3$, $\Gamma(\mathbb{Z}_{p}[x]/\langle x^{4} \rangle)$ is hypoerenergetic.
		\end{itemize}
\end{proposition}

\begin{figure}[H]
	\begin{center}
			\begin{minipage}{.4\textwidth}
		\centering
		\scalebox{.3}{\includegraphics{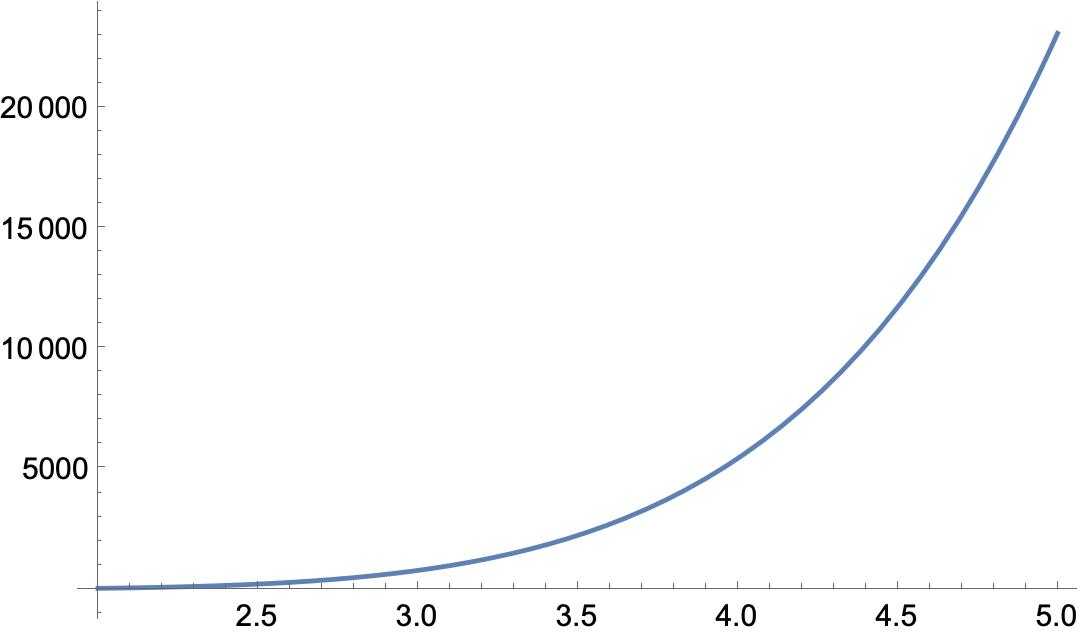}} 
		$f(p)$
	\end{minipage}    
	\begin{minipage}{.3\textwidth}
		\centering
		\scalebox{.3}{\includegraphics{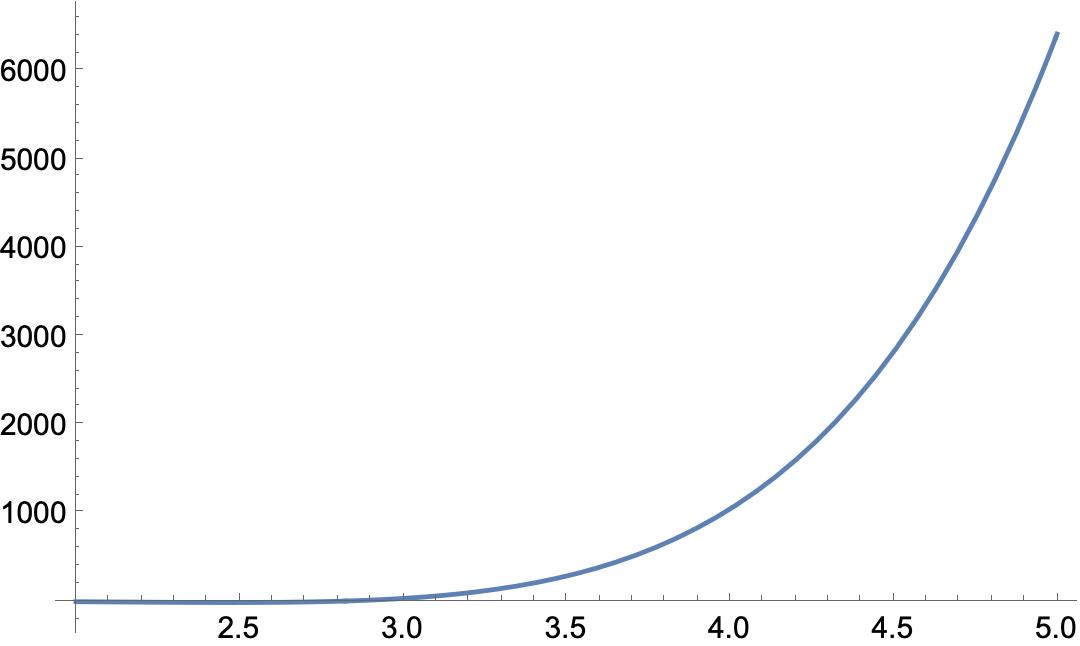}}
		$g(p)$
	\end{minipage}
	\end{center}
	\caption{Picture showing the increasing behaviour of $f(p)$ and $g(p)$. }
	\label{check}
\end{figure}

\noindent\textbf{Proof.} Let $\Gamma(\mathbb{Z}_{p}[x]/\langle x^{4} \rangle) $ be the zero divisor graph of order $n=p^{3}-1.$ Also, it is well known that the energy of the complete graph of order $n=p^{3}-1$ is $\mathcal{E}(K_{n})=2(p^{3}-2).$ By Theorem \ref{energy theorem}, the energy of $\Gamma(\mathbb{Z}_{p}[x]/\langle x^{4} \rangle)$ satisfies
\[ \mathcal{E}(\Gamma(\mathbb{Z}_{p}[x]/\langle x^{4} \rangle))\leq p^{2}-3+\sqrt{9p^{4}-12p^{3}-6p^{2}+15}. \]
Suppose to the contrary $ \Gamma(\mathbb{Z}_{p}[x]/\langle x^{4} \rangle)$ is hyperenergetic, that is,
\[ \mathcal{E}(\Gamma(\mathbb{Z}_{p}[x]/\langle x^{4} \rangle))> 2(p^{3}-2), \]
that is, \[ p^{2}-3+\sqrt{9p^{4}-12p^{3}-6p^{2}+15}>2(p^{3}-2), \]
which after simplification gives 
\begin{equation}\label{eq check}
	2 p^6-2 p^5-4 p^4+4 p^3+4 p^2-7<0.
\end{equation}
Clearly Inequality \eqref{eq check} is not true, since $ f(p)=2 p^6-2 p^5-4 p^4+4 p^3+4 p^2-7$ is an increasing function for primes $p\geq 2.$ In fact the global minimum attained by $f(p)$ is $41$ for $p\geq 2$.  The graphical representation of Inequality \eqref{eq check} is shown in Figure \ref{check}
Therefore, $\mathcal{E}(\Gamma(\mathbb{Z}_{p}[x]/\langle x^{4} \rangle))< 2p^{3}-2$ and hence is non-hyperenergetic for any prime $p.$\vskip 2mm
(b). Comparing $\mathcal{E}(\Gamma(\mathbb{Z}_{p}[x]/\langle x^{4} \rangle))\leq  p^{2}-3+\sqrt{9p^{4}-12p^{3}-6p^{2}+15}< p^{3}-1,$ we obtain 
\[ g(p)=p^6-2p^{5}-8p^{4}+16p^{3}+2p^{2}-11>0, \]
which is always true for any prime $p\geq 3$ and is negative for $p\in [2,3).$
The plot of $g(p)$ is shown in Figure \ref{check}.
\qed
\section{Topological indices of zero divisor graphs of $ \Gamma(\mathbb{Z}_{p}[x]/\langle x^{4} \rangle)$}\label{section 3}
\paragraph{}
For a graph $ G, $ a general vertex-degree-based (VDB) topological index $ \phi $ \cite{gutman2022} is defined as
\[ \phi(G)=\sum_{uv\in E(G)}\phi_{d_{u}d_{v}}, \]
$ \phi_{d_{u}d_{v}} $ is a function with the property $ \phi_{d_{u}d_{v}}=\phi_{d_{v}d_{u}} $. For particular values of $ \phi_{d_{u}d_{v}} $, we get well know topological indices like arithmetic-geometric index \cite{shegehalli} $ \phi_{d_{u}d_{v}}=\frac{d_{u}+d_{v}}{2\sqrt{d_{u}d_{v}}} $, the general Randi\'c index \cite{randic1} $ \phi_{d_{u}d_{v}}=(d_{u}d_{v})^{\alpha}, $ $\Big($for $ \alpha =-\frac{1}{2} $, we obtain the ordinary Randi\'c index \cite{randic2} $R=\displaystyle\sum_{uv\in E(G)} \frac{1}{\sqrt{d_{u}d_{v}}},$ and for $\alpha=1$ we get the second Zagreb index \cite{gutmanzagreb}, $M_{2}(G)=\sum_{uv\in E(G)}(d_{u}\cdot d_{v})$  $\Big)$, the general Sombor index $ \phi_{d_{u}d_{v}}= \left(d_{u}^{2}+d_{v}^{2} \right)^{\alpha}, $ the general sum connectivity index \cite{zhou} with $ \phi_{d_{u}d_{v}}=(d_{u}+d_{v})^{\alpha}$ and for $\alpha=1,$ we get the first Zagreb index \cite{gutmanzagreb}, $M_{1}(G)=\sum_{uv\in E(G)}(d_{u}+d_{v})$, the general inverse sum indeg index \cite{gutman2109} $\phi_{d_{u}d_{v}}=\left( \frac{d_{u}d_{v}}{d_{u}+d_{v}}\right)^{\alpha} $ and for $\alpha=1,$ we get the inverse inverse sum indeg index.
 Similarly, all the degree based topological indices can be generated by the function $ \phi_{d_{u}d_{v}}.$ The exponential of the topological indices as given as $ e^{\phi}(G)=\sum_{uv\in E(G)}e^{\phi_{d_{u}d_{u}}} $ \cite{rada}.
 A topological index of $G$ is an entity (numeric) related to a graph structure that links chemical structure with certain physico-chemical relation, like biological activity or chemical reactivity.  With molecular modelling, the relationship between the structural properties and activity of chemical compounds can be recorded/observed \cite{gutmanacta}. The aim of descriptors were discovered to play a significant role in (QSPR/QSAR) analysis of molecular structures, earlier Wiener in 1947 \cite{wiener}.

The following result gives a general formulae for the topological indices for the zero divisor graph of $\mathbb{Z}_{p}[x]/\langle x^{4} \rangle.$
\begin{theorem}\label{vdb topological indices}
	For $G\cong \Gamma(\mathbb{Z}_{p}[x]/\langle x^{4} \rangle),$ the VDB topological index is
	\[\binom{p-1}{2}\phi_{d_{u}d_{u}}+(p-1)(p^{2}-p)\phi_{d_{u}d_{v}}+\binom{p^{2}-p}{2}\phi_{d_{v}d_{v}}+(p^{3}-p^{2})(p-1)\phi_{d_{u}d_{w}}, \]
	where $d_{u}$ represent the common degree of the vertices in $A,$ $d_{v}$ represent the common degree of the vertices in $B$ and $d_{w}$ represent the common degree of the vertices in $C.$
\end{theorem}
\noindent\textbf{Proof.} Since the vertices of $G$ are partitioned into three mutually disjoint sets $A,$ $B$ and $C.$ The vertices in $A$ form a clique of order $p-1$, the vertices in $B$ form a clique of size $p^{2}-p$ and each vertex of $A$ is connected to every vertex of $B.$ Also, $C$ is an independent set of cardinality $p^{3}-p$ and each such vertex is adjacent to every vertex of $ A.$ With this structure, the number of edges inside $A$ are $\binom{p-1}{2},$ that of in $B$ are $\binom{p^{2}-p}{2},$ between $A$ and $B$, there are $(p-1)(p^{2}-p)$ edges and between $A$ and $C$, there are $(p-1)(p^{3}-p^{2})$ edges. Let the vertex labelling of $G$ be
\[ \big\{u_{1}, u_{2}, \dots, u_{|A|}, v_{1}, v_{2}, \dots, v_{|B|}, w_{1}, w_{2},\dots, w_{|C|} \big\}, \]
where $u_{i}$'s are the vertices of $A$, $v_{i}$'s are the vertices of $B$ and $w_{i}$'s are the vertices of $C.$ From the structure of $G,$ it is clear that the vertices in $A$ have common degree $p^{3}-2$ and we denote it by  $d_{u}=p^{3}-2.$ Similarly, the vertices in $B$ and $C$ have common degrees and are denoted by $d_{v}=p^{2}-2$ and $d_{w}=p-1$, respectively. Thus, the VDB topological index of $G$ is
\begin{align*}
	\phi(G)&=\sum_{u_{i}u_{j}\in E(G)}\phi_{d_{u_{i}}d_{u_{j}}}+\sum_{u_{i}v_{j}\in E(G)}\phi_{d_{u_{i}}d_{v_{j}}}+\sum_{v_{i}v_{j}\in E(G)}\phi_{d_{u_{i}}d_{u_{j}}}+\sum_{u_{i}w_{j}\in E(G)}\phi_{d_{u_{i}}d_{w_{j}}}\\
	&=\binom{p-1}{2}\phi_{d_{u}d_{u}}+(p-1)(p^{2}-p)\phi_{d_{u}d_{v}}+\binom{p^{2}-p}{2}\phi_{d_{v}d_{v}}+(p^{3}-p^{2})(p-1)\phi_{d_{u}d_{w}}.
\end{align*}
\qed
The following is an immediate consequence of Theorem \ref{vdb topological indices}, which gives the closed formulae for the first and the second Zagreb index of $\Gamma(\mathbb{Z}_{p}[x]/\langle x^{4} \rangle)$.
\begin{corollary}\label{vdb corollary}
	For $G\cong \Gamma(\mathbb{Z}_{p}[x]/\langle x^{4} \rangle),$ the following hold.
	\begin{itemize}
		\item[\bf (i)] $M_{1}(G)=p^{7}-11p^{4}+11p^{3}+3p^{2}-4.$
		\item[\bf (ii)] $M_{2}(G)=\frac{1}{2}\Big( 6p^{8} - 15p^{7} + 2p^{6} + 3p^{5} + 32p^{4} - 28p^{3} - 8p^{2} + 8\Big).$
	\end{itemize}
\end{corollary}
\noindent\textbf{Proof.} By substituting $\phi_{d_{u}d_{v}}=d_{u}+d_{v}$ and $ \phi_{d_{u}d_{v}}=d_{u}\cdot d_{v}$ in Theorem \ref{vdb topological indices}, respectively. The result follows. \qed

We state the result of \cite{johnson} for the Zagreb indices $M_{1}$ and $M_{2}$ for the zero divisors of $\Gamma(\mathbb{Z}_{p}[x]/\langle x^{4} \rangle). $
\begin{theorem}[Theorem 3.4, \cite{johnson}]\label{johnsom theorem zagreb}
	Let $R = \mathbb{Z}_{p}[x]/\langle x^{4}\rangle$ with prime $p$, then
\begin{align*}
		M_{1}(\Gamma(\mathbb{Z}_{p}[x]/\langle x^{4}\rangle))&= p^{7}-11p^{4} + 11p^{3} + 3p^{2}- 4\\
	M_{2}(\Gamma(\mathbb{Z}_{p}[x]/\langle x^{4}\rangle))&= 12(6p^{8}-15p^{7} + 2p^{6} + 3p^{5} + 32p^{4} + 28p^{3}-8p^{2}-8).
\end{align*}
\end{theorem}

As we observe that the values for $M_{2}(G)$ are different in Theorem \ref{johnsom theorem zagreb} than given in Corollary \ref{vdb corollary}. We observe that there is calculation correction in Theorem 3.4 \cite{johnson}, in fact in the proof of Theorem \ref{johnsom theorem zagreb}, the evaluation of $M_{2}(G)$ in third line is written incorrectly. We verify the second Zagreb index by software and by Corollary \ref{vdb corollary} for the zero divisor graph $ \Gamma(\mathbb{Z}_{3}[x]/\langle x^{4}\rangle)$ in the following remark.
\begin{remark}
	For the graph $G\cong \Gamma(\mathbb{Z}_{3}[x]/\langle x^{4} \rangle) $ shown in Figure \ref{fig 1}. We calculate the values of $M_{1}$ and $M_{2}.$ With $p=3$ in Corollary \ref{vdb corollary}, the first and seconds Zagreb indices are
	\[ M_{1}(G)=1616~\text{and}~ M_{2}(G)=5260. \]
	While by  By Theorem \ref{johnsom theorem zagreb} (Theorem 3.4, \cite{johnson}) with $p=3$, the value of the second Zagreb index is $M_{2}(G)=6008,$ but by software calculation (Mathematica, SageMath or AutographiX) and by Corollary, the correct value of $M_{2}(G) $ is $5260.$ Thus, the formulae for $M_{2}(\Gamma(\mathbb{Z}_{p}[x]/\langle x^{4}\rangle))$ given in Theorem \ref{johnsom theorem zagreb} (Theorem 3.4, \cite{johnson}) is not correct.
\end{remark}

Hansen and Vuki$\check{c}$cevi\'c \cite{hansen} compared the  first and the second Zagreb indices, and  put forward the following conjecture.
\begin{conjecture}\label{conjecture}
	For a simple graph $G,$
	\begin{equation}\label{conjecture eq}
		 \frac{M_{2}(G)}{|E(G)|}\geq \frac{M_{1}(G)}{|V(G)|}.
	\end{equation}
\end{conjecture}

The authors in \cite{hansen} showed that Conjecture \ref{conjecture} holds for chemical graphs. In \cite{vuk}, it was shown that the conjecture holds for trees with equality for a star graph. In \cite{liu}, the authors settled the conjecture \ref{conjecture}  for connected unicyclic graphs with equality when the graph is a cycle. The  equality case of \ref{conjecture eq} is studied in \cite{vuk1}. A survey on comparing Zagreb indices can be seen in \cite{liu1}.

For the graph $G\cong \Gamma(\mathbb{Z}_{p}[x]/\langle x^{4}\rangle)$. The second Zagreb Index given in Corollary \ref{vdb corollary} can be written as
\[ M_{2}(G)=\frac{1}{2} (p-1) \left(6 p^7-9 p^6-7 p^5-4 p^4+28 p^3-8 p-8\right), \]
the edge cardinality of $G$ is
\[ |E(G)|=\frac{1}{2} (p-1) \left(3 p^3-p^2-2 p-2\right). \]
By Inequality \ref{conjecture eq}, we have  
\[ \dfrac{6 p^7-9 p^6-7 p^5-4 p^4+28 p^3-8 p-8}{3 p^3-p^2-2 p-2}\geq \frac{p^{7}-11p^{4}+11p^{3}+3p^{2}-4}{p^{3}-1}. \]
This gives us
\begin{equation}\label{conj check 1}
	 (p-1) p^2 \left(3 p^7-5 p^6-10 p^5+15 p^4+8 p^3-5 p^2-6 p-2\right)\geq 0.
\end{equation}
Since $p-1$ and $p^{2}$ are positive for $p\geq 2,$ so we consider $$h(p)=3 p^7-5 p^6-10 p^5+15 p^4+8 p^3-5 p^2-6 p-2.$$ The global minimum value attained by $h(p)$ is $14$ at $p=2,$ and the function is strictly increasing for $p\geq 3.$ Thus, Inequality \ref{conj check 1} is true and it follows that Inequality \eqref{conjecture eq} is true for $G\cong \Gamma(\mathbb{Z}_{p}[x]/\langle x^{4}\rangle).$ Hence, we have the following result.
\begin{proposition}
	Let $G\cong \Gamma(\mathbb{Z}_{p}[x]/\langle x^{4}\rangle)$ be the zero divisor graph. Then the first Zagreb index $M_{1}(G)$ and the second Zagreb index $M_{2}(G)$  satisfy Conjecture \ref{conjecture}.
\end{proposition}

The last result is a consequence of Theorem \ref{vdb topological indices} and gives some well known general topological indices of $\Gamma(\mathbb{Z}_{p}[x]/\langle x^{4}\rangle).$
\begin{corollary}\label{vdb corollary 2}
	For $G\cong \Gamma(\mathbb{Z}_{p}[x]/\langle x^{4} \rangle),$ with $d_{u}=p^{3}-2, d_{v}=p^{2}-2$ and $d_{w}=p-1,$
	 the following hold.
	\begin{itemize}
		\item[\bf (i)] The general Randi\'c index  of $G$ is $R_{\alpha}(G)=\binom{p-1}{2}(d_{u}^{2})^{\alpha}+\binom{p^{2}-p}{2}(d_{v}^{2})^{\alpha}+p(p-1)^{2}(d_{u}d_{v})^{\alpha}+p^{2}(p-1)^{2}(d_{v}d_{w})^{\alpha}.$
		\item[\bf (ii)] The general sum-connectivity index  of $G$ is $\chi_{\alpha}(G)=\binom{p-1}{2}(2d_{u})^{\alpha}+\binom{p^{2}-p}{2}(2d_{v})^{\alpha}+p(p-1)^{2}(d_{u}+d_{v})^{\alpha}+p^{2}(p-1)^{2}(d_{v}+d_{w})^{\alpha}.$
		\item[\bf (iii)] The general inverse sum indeg index  of $G$ is $S_{\alpha}(G)=\binom{p-1}{2}\left (\frac{d_{u}}{2}\right )^{\alpha}+\binom{p^{2}-p}{2}\left (\frac{d_{v}}{2}\right )^{\alpha}+p(p-1)^{2}\left (\frac{d_{u}d_{v}}{d_{u}+d_{v}}\right )^{\alpha}+p^{2}(p-1)^{2}\left (\frac{d_{v}d_{w}}{d_{v}+d_{2}}\right )^{\alpha}.$
		\item[\bf (iv)] The general Sombor index  of $G$ is $SO_{\alpha}(G)=\binom{p-1}{2}(2d_{u}^{2})^{\alpha}+\binom{p^{2}-p}{2}(2d_{v}^{2})^{\alpha}+p(p-1)^{2}(d_{u}^{2}+d_{v}^{2})^{\alpha}+p^{2}(p-1)^{2}(d_{v}^{2}+d_{w}^{2})^{\alpha}.$
	\end{itemize}
\end{corollary}
\noindent\textbf{Proof.} The proof follows by taking $\phi_{d_{u}d_{u}}=(d_{u}d_{v})^{\alpha}, \phi_{d_{u}d_{u}}=(d_{u}+d_{v})^{\alpha}, \phi_{d_{u}d_{u}}=\left (\frac{d_{u}d_{v}}{d_{u}+d_{v}}\right )^{\alpha}, $ and $ \phi_{d_{u}d_{u}}=(d_{u}^{2}+d_{v}^{})^{\alpha}$ in Theorem \ref{vdb topological indices}.

 \section{Conclusion} \label{section 4}
 \paragraph{}
 The present study corrects the results related to the characteristic polynomial, the energy and the Zagreb indices of the zero divisor graphs of $\mathbb{Z}_{p}[x]/\langle x^{4} \rangle,$ thereby correcting the main results of Johnson and  Sankar \cite{johnson}. We also prove that Hansen and Vuki$\check{c}$cevi\'c conjecture holds for $\Gamma(\mathbb{Z}_{p}[x]/\langle x^{4} \rangle)$ and present a closed formulae so that all the well known topological indices of $\Gamma(\mathbb{Z}_{p}[x]/\langle x^{4} \rangle)$ can be derived from aforesaid result.

\vskip 3mm
\noindent\textbf{Data Availability:}\\
There is no data associated with this article.
\vskip 3mm

\noindent\textbf{Conflict of interest}\\
The authors declare that they have no competing interests.

\end{document}